\documentclass[letterpaper, 11pt,  reqno]{amsart}

\usepackage{amsmath,amssymb,amscd,amsthm,amsxtra, esint}

\usepackage[english, french]{babel}

\usepackage{hyperref}

\usepackage{cases}

\allowdisplaybreaks[2]

\usepackage{color}

\definecolor{gr}{rgb}   {0.,   0.69,   0.23 }
\definecolor{bl}{rgb}   {0.,   0.5,   1. }
\definecolor{mg}{rgb}   {0.85,  0.,    0.85}
\definecolor{yl}{rgb}   {0.8,  0.7,   0.}
\definecolor{or}{rgb}  {0.7,0.2,0.2}

\catcode`\Ž = \active\def Ž{\'e} \catcode`\ƒ = \active\def ƒ{\'E}  \catcode`\Ë = \active\def Ë{\`A}\catcode`\ = \active\def {\`e}
\catcode`\ = \active\def {\c{c}} \catcode`\ˆ = \active\def
ˆ{\`a} \catcode`\ = \active \def {\^e} \catcode`\ = \active
\def {\`u} \catcode`\™ = \active \def ™{\^o} \catcode`\ž =
\active \def ž{\^u} \catcode`\" = \active \def "{\^{\i}}
\catcode`\‰ = \active \def ‰{\^a} \catcode`\š = \active \def
š{{\"o}}

\newtheorem{theorem}{Theorem} [section]

\newtheorem{lemma}[theorem]{Lemma}
\newtheorem{proposition}[theorem]{Proposition}
\newtheorem{remark}[theorem]{Remark}

\newcommand{\noi}{\noindent}
\newcommand{\Z}{\mathbb{Z}}
\newcommand{\R}{\mathbb{R}}
\newcommand{\C}{\mathbb{C}}
\newcommand{\T}{\mathbb{T}}

\let\P= \undefined
\newcommand{\P}{\mathbf{P}}

\newcommand{\E}{\mathbb{E}}

\newcommand{\F}{\mathcal{F}}

\newcommand{\be}{\beta}
\newcommand{\dl}{\delta}

\newcommand{\nb}{\nabla}

\newcommand{\Dl}{\Delta}
\newcommand{\eps}{\varepsilon}

\newcommand{\g}{\gamma}
\newcommand{\G}{\Gamma}
\newcommand{\ld}{\lambda}

\newcommand{\s}{\sigma}
\newcommand{\Si}{\Sigma}
\newcommand{\ft}{\widehat}

\newcommand{\wt}{\widetilde}
\newcommand{\cj}{\overline}
\newcommand{\dx}{\partial_x}
\newcommand{\dt}{\partial_t}
\newcommand{\dd}{\partial}

\renewcommand{\l}{\ell}
\renewcommand{\o}{\omega}
\renewcommand{\O}{\Omega}

\newcommand{\les}{\lesssim}
\newcommand{\ges}{\gtrsim}

\newcommand{\jb}[1]
{\langle #1 \rangle}

\newcommand{\ind}{\mathbf 1}

\newcommand{\M}{\mathcal{M}}

\newcommand{\N}{\mathbb{N}}

\newcommand{\Pk}{P^{(2m+2)}_2}
\newcommand{\Pkn}{P^{(2m+2)}_{2, N}}
\renewcommand{\H}{\mathcal{H}}

\newtheorem*{ackno}{Acknowledgements}

\numberwithin{equation}{section}
\numberwithin{theorem}{section}

\begin{document}
\baselineskip = 15pt

\selectlanguage{english}

\title[Invariant Gibbs measures for the 2-$d$ defocusing NLW]
{Invariant Gibbs measures for 
the 2-$d$ defocusing  nonlinear wave equations}

\author[T.~Oh and  L.~Thomann]
{Tadahiro Oh and Laurent Thomann}

\address{
Tadahiro Oh, School of Mathematics\\
The University of Edinburgh\\
and The Maxwell Institute for the Mathematical Sciences\\
James Clerk Maxwell Building\\
The King's Buildings\\
Peter Guthrie Tait Road\\
Edinburgh\\ 
EH9 3FD\\
 United Kingdom}

\email{hiro.oh@ed.ac.uk}

\address{
Laurent Thomann\\
Institut  \'Elie Cartan, Universit\'e de Lorraine, B.P. 70239,
F-54506 Vand\oe uvre-l\`es-Nancy Cedex, France}

\email{laurent.thomann@univ-lorraine.fr}



\maketitle

\vspace{-8mm}

\begin{abstract}
We consider the defocusing nonlinear wave equations (NLW)
on the two-dimensional torus.
In particular, we construct invariant Gibbs measures
for the renormalized so-called Wick ordered NLW.
We then prove weak universality of the Wick ordered NLW,
showing that the Wick ordered NLW naturally appears as a suitable scaling limit
of non-renormalized NLW with Gaussian random initial data.

\end{abstract}

\vspace{-3mm}

\begin{otherlanguage}{french}
\begin{abstract}
On consid\`ere les \'equations des ondes non-lin\'eaires d\'efocalisantes sur le tore de dimension deux. On construit des mesures de Gibbs invariantes pour  les Žquations renormalisŽes au sens de Wick. 
On prouve ensuite une propri\'et\'e d'universalit\'e  faible pour ces \'equations renormalisŽes, en montrant qu'elle apparaissent comme limites d'Žquations d'ondes non renormalis\'ees avec conditions initiales al\'eatoires de loi gaussienne.

\end{abstract}
\end{otherlanguage}
%

%
%

\baselineskip = 14pt

\section{Introduction}

\subsection{Nonlinear wave equations}

We consider 
the defocusing nonlinear wave equations
 (NLW) in two spatial dimensions:
\begin{align}
\begin{cases}
\dt^2 u -  \Dl  u + \rho u  +  u^{2m+1} = 0 \\
(u, \dt u) |_{t = 0} = (\phi_0, \phi_1), 
\end{cases}
\qquad (t,x) \in \R \times \M,
\label{NLW1}
\end{align}

\noi
where $\rho \geq 0$ and  $m \in \N $.
When $\rho > 0$, \eqref{NLW1}
is also referred to as the nonlinear Klein-Gordon equation.
We, however, simply refer to \eqref{NLW1} as NLW
and moreover restrict our attention to the real-valued setting.
In the following, we mainly consider \eqref{NLW1}
on the two-dimensional torus $\M = \T^2 = (\R/\Z)^2$
but we also provide a brief discussion  when $\M$ 
is  a two-dimensional compact Riemannian manifold without boundary or a bounded domain in~$\R^2$
(with the Dirichlet or Neumann boundary condition).
See Theorem~\ref{THM:GWP2} below.

Our main goal in this paper is 
to construct  an invariant  Gibbs measure for a renormalized version of~\eqref{NLW1}
by studying dynamical properties of the renormalized equation.

\subsection{Gibbs measures and Wick renormalization}
With $v = \dt u$, 
we can write the equation~\eqref{NLW1} 
in the following Hamiltonian formulation:
\begin{equation*}
 \dt 
 \begin{pmatrix}
 u \\ v 
 \end{pmatrix}
 =  
 \begin{pmatrix}
 0& 1 \\ -1 & 0
 \end{pmatrix}
\frac{\dd H}{\dd(u, v )}, 
\end{equation*}

\noi
where $H = H(u, v)$ is  the Hamiltonian given by 
\begin{align}
H(u, v) = \frac{1}{2}\int_{\T^2}\,\big(\rho u^2 +  |\nb u|^2\big)\, dx
+ 
\frac{1}{2}\int_{\T^2} v^2dx
+ \frac1{2m + 2} \int_{\T^2} u^{2m+2} dx.
\label{Hamil}
\end{align}

\noi
By drawing an analogy to the finite dimensional setting,
the Hamiltonian structure of the equation and the conservation of the Hamiltonian
suggest that the Gibbs measure $\Pk$ of the form:\footnote{Henceforth, we use  
$Z$,  $Z_N$, etc.~to denote various  normalizing constants
so that the corresponding measures are probability measures when appropriate.}
\begin{align}
 \text{``}d\Pk = Z^{-1} \exp(- \be H(u,v ))du\otimes dv \text{''}
\label{G1}
 \end{align}

\noi
is  invariant under the dynamics of~\eqref{NLW1}.\footnote{
We simply set $\be = 1$ in the following.
While our analysis holds for any $\be > 0$, 
the resulting (renormalized) Gibbs measures are mutually singular
for  different values of $\be > 0$.
See~\cite{OQ}.}
With \eqref{Hamil}, we can rewrite the formal expression~\eqref{G1}
as 
\begin{align}
d\Pk
& = Z^{-1} 
e^{-\frac 1{2m+2} \int u^{2m+2} dx} 
e^{-\frac 12 \int (\rho u^2 + |\nb u |^2) dx} du\otimes e^{-\frac 12 \int v^2 dx} dv \notag\\
& \sim e^{-\frac 1{2m+2} \int u^{2m+2} dx} d\mu , 
\label{G2}
 \end{align}

\noi
where $\mu$ is  the Gaussian measure $\mu$ on $\mathcal{D}'(\T^2)\times \mathcal{D}'(\T^2)$ 
with the density\footnote{On $\T^2$, we need to  assume $\rho > 0$ in order to avoid a problem at the zeroth frequency. See \eqref{G5} below.
In the case of a bounded domain in $\R^2$ with the Dirichlet boundary condition, 
we can take $\rho = 0$.}
\begin{align}
 d\mu =  Z^{-1} e^{-\frac 12 \int( \rho u^2 + |\nb u |^2) dx} du\otimes e^{-\frac 12 \int v^2dx} dv.
\label{G3}
\end{align}

\noi
Note that $\mu$ has a tensorial structure: $\mu = \mu_0 \otimes \mu_1$, where
the marginal measures~$\mu_0$ and~$\mu_1$ are given by 
\begin{align}
 d\mu_0 =  Z_0^{-1} e^{-\frac 12 \int (\rho u^2 + |\nb u |^2 )dx} du
 \qquad \text{and}\qquad 
 d\mu_1 =  Z_1^{-1} e^{-\frac 12 \int v^2dx} dv.
\label{G4}
\end{align}

\noi
Namely, $\mu_0$ is the Ornstein-Uhlenbeck measure
and $\mu_1$ is the white noise measure on $\T^2$.

Recall that $\mu$ is the induced probability measure
under the map:\footnote{We drop the harmless factor $2\pi$ in the following.}
\begin{align}
 \o \in \O \longmapsto (u, v) 
 = \bigg( \sum_{n \in \Z^2} \frac{g_{0, n}(\o)}{\jb{n}_\rho}e^{in\cdot x},
\sum_{n \in \Z^2} g_{1, n}(\o)e^{in\cdot x}\bigg), 
\label{G5}
 \end{align}

\noi
where 
$\jb{n}_\rho = \sqrt{\rho+|n|^2}$ and 
$\{g_{0, n}, g_{1, n}\}_{n \in \Z^2}$
is a sequence of independent standard complex-valued Gaussian
random variables on a probability space $(\O, \F, P)$
conditioned that ${g_{j, -n} = \cj{g_{j, n}}}$, $n \in \Z^2$, $j = 0, 1$.
In view of \eqref{G5}, it is easy to see that $\mu$ is supported
on \[\H^s(\T^2) : = H^s(\T^2)\times H^{s-1}(\T^2) , \quad s < 0.\]

\noi
Moreover, we have $\mu(\H^0(\T^2)) = 0$.
This  implies that 
$\int u^{2m+2} dx = \infty$
almost surely with respect to $\mu$.
In particular,  the right-hand side of \eqref{G2} would not be a probability measure, 
thus requiring a {\it renormalization} of the potential part of the Hamiltonian.
In the two-dimensional case, it is known that 
a Wick ordering suffices for this purpose.
See Simon~\cite{Simon} and  Glimm-Jaffe~\cite{GJ}.
Also, see Da Prato-Tubaro~\cite{DPT1}
for a concise discussion on $\T^2$,
where the Gibbs measures naturally appear in the context of the stochastic
quantization equation.

In the following, we give a brief review of the Wick renormalization on $\T^2$.
See~\cite{DPT1} for more details. 
Let $u$ denote a typical element under  $\mu_0$ defined in~\eqref{G4}.
Since $u \notin L^2(\T^2)$ almost surely, we have 
\[ \int_{\T^2} u^2 dx = \lim_{N \to \infty} \int_{\T^2} (\P_N u)^2 dx = \infty\]

\noi
almost surely,
where 
 $\P_N$ is the Dirichlet projection onto the frequencies $\{|n|\leq N\}$.

For each $x \in \T^2$,  
 $\P_N u(x)$ 
 is a mean-zero real-valued Gaussian random variable with variance\footnote{Note that $\s_N$ 
 defined in \eqref{G6} is independent of $x \in \T^2$.
When $\M$ is a two-dimensional compact Riemannian manifold without boundary or a bounded domain in $\R^2$,
the variance  $\s_N(x) = \E [(\P_N u)^2(x)]$  depends on $x \in \M$
but satisfies the logarithmic bound in $N$.  See \eqref{4sigma} below.
}
\begin{align}
\s_N \stackrel{\text{def}}{=} \E [ (\P_N u)^2(x)]
= \sum_{|n|\leq N}\frac{1}{\rho+|n|^2} \sim \log N.
\label{G6}
\end{align}

\noi
This motivates us to define the Wick ordered monomial
 $:\! (\P_N u)^k  \!:$
by 
\begin{align}
:\! (\P_N u)^k (x)\!:  \, = H_k(\P_N u(x); \s_N)
\label{Wick1}
\end{align}

\noi
in a pointwise manner.
Here,  $H_k(x; \s)$ is the Hermite polynomial of degree $k$
defined in~\eqref{H1}.
Then, with \eqref{G5} and \eqref{G6},  it is easy to see that 
the random variables $X_N(u)$ defined by 
\begin{align*}
X_N(u)  = \int_{\T^2} 
:\!  (\P_N u)^2 (x)\!: dx
\end{align*}

\noi
have uniformly bounded second moments
and converge to some random variable in $L^2(d\mu_0)$ which we denote by 
\begin{align*}
X_\infty(u)  = \int_{\T^2} :\! u^2 \!: dx \in L^2(d\mu_0).   
\end{align*}

\noi
In view of the Wiener chaos estimate (Lemma \ref{LEM:hyp3}), we see that  $X_N(u)$ also converges to $X_\infty(u)$ in 
$L^p(d\mu_0)$, $p < \infty$.

In general, given any $m \in \N$, 
one can show that the limit
\begin{align}
\int_{\T^2} :\! u^{2m+2}\!: dx
= \lim_{N \to \infty} 
\int_{\T^2}:\! (\P_N u)^{2m+2}\!: dx
\label{Wick2}
\end{align}
	
\noi
exists in $L^p( \mu)$
for any finite $p \geq 1$.
Moreover, we have the following proposition.

\begin{proposition}\label{PROP:Gibbs1}
Let $m \in \N$.
Then, 
$R_N(u)\stackrel{\textup{def}}{=}e^{-\frac 1{2m+2} \int_{\T^2} :(\P_N u)^{2m+2} : \, dx} \in L^p(\mu)$ 
for any finite $p\geq 1$
with a uniform bound in $N$, depending on $p \geq 1$.
Moreover, for any finite $p \geq 1$, 
$R_N(u)$ converges to some $R(u)$ in $L^p(\mu)$ as $N \to \infty$.
\end{proposition}

This proposition follows from 
the hypercontractivity of the Ornstein-Uhlenbeck semigroup
and
Nelson's estimate~\cite{Nelson2}. 
See also \cite{DPT1, OT}.
Denoting the limit $R(u)\in L^p(\mu)$ by
\[ R(u) =  e^{-\frac 1{2m+2} \int_{\T^2} : u^{2m+2} : \, dx},\]

\noi
Proposition~\ref{PROP:Gibbs1}
allows us to define the Gibbs measure $\Pk$  
associated with 
the Wick  ordered Hamiltonian:
\begin{align*}
H_\text{\tiny Wick}(u, v) =\frac{1}{2}\int_{\T^2}\,\big(\rho u^2 +  |\nb u|^2\big)\, dx
+ 
\frac{1}{2}\int_{\T^2} v^2dx
+ \frac1{2m+2} \int_{\T^2} :\! u^{2m+2}\!: dx
\end{align*}

\noi
by 
\begin{align*}
d \Pk 
& = Z^{-1}e^{-H_{\text{Wick}}(u, v)} du\otimes dv
=  Z^{-1}  e^{-\frac 1{2m+2} \int_{\T^2} :u^{2m+2} : \,dx} d\mu\notag\\
& =Z^{-1} R(u) d\mu.
 \end{align*}

\noi
It follows from Proposition~ \ref{PROP:Gibbs1}
that   $\Pk \ll \mu$  and, in particular, $\Pk$ is a probability measure on $\H^s(\T^2)\setminus \H^0(\T^2)$, $s< 0$.
Moreover, 
defining $\Pkn$ 
by \[d \Pkn = Z_N^{-1} R_N(u) d\mu ,\]

\noi 
we see that 
$\Pkn$ converges ``uniformly'' to $\Pk$
in the sense that 
given any  $\eps > 0$, there exists $N_0 \in \N$ such that 
\[ \big|  \Pkn(A) - \Pk(A)\big| < \eps \]

\noi
for any $N \geq N_0$
and  any measurable set 
$A \subset \H^s(\T)$, $s<0$.

Lastly, let us briefly discuss 
the construction of  the Gibbs measure $\Pk$
when $\M$ is  a two-dimensional compact Riemannian manifold without boundary or a bounded domain in $\R^2$
(with the Dirichlet or Neumann boundary condition).
In this case, the Gaussian measure $\mu$ in \eqref{G3} 
represents the induced probability measure under the map:
\begin{align}
 \o \in \O \longmapsto 
 (u, v)  = \bigg( \sum_{n \in \N }\frac{g_{0, n}(\omega)}{(\rho+\lambda^2_n)^{\frac{1}{2}}}\, \varphi_n(x),
  \sum_{n \in \N }g_{1, n}(\omega) \varphi_n(x)\bigg),
\label{4G5}
 \end{align}

\noi
where 
 $\{\varphi_n\}_{n \in \N }$  is an orthonormal basis   of $L^2(\M)$ consisting  of eigenfunctions of 
 the Laplace-Beltrami operator $-\Delta$
with 
the corresponding eigenvalues $\{\lambda_n^2\}_{n \in \N }$, which we assume to be arranged in the increasing order.
It is easy to see from \eqref{4G5} that $\mu$ is supported on 
$\H^s(\M)\setminus \H^0(\M)$, $s <0$.

Given $N\in \N$, 
we define $\s_N$   by
\begin{equation}
\s_N (x)  = \E [(\P_N u_N)^2(x)]
= \sum_{\lambda_{n}\leq N}\frac{\varphi_n(x)^2 }{\rho+\lambda_{n}^{2}} \les \log N,
\label{4sigma}
\end{equation}

\noi
where $\P_{N}$ denotes the spectral projector defined by 
\[\P_N u= \sum_{\lambda_{n}\leq N }\ft u(n) \varphi_{n}.\]

\noi
Note that unlike the situation on $\T^2$, 
$\s_N(x)$ now depends on $x \in \M$.
The last inequality in \eqref{4sigma}, however, holds
independently of $x \in \M$
thanks to 
 Weyl's law  
$\lambda_n\approx n^{\frac 12}$ 
(see \cite[Chapter~14]{Zworski})
and \cite[Proposition 8.1]{BTT1}.
With this definition of $\s_N(x)$, we can define the
Wick ordered monomials  $:\! (\P_N u)^k  \!:$
as in \eqref{Wick1}
and  $:\! u^k  \!:$ by the limiting procedure.
Then, the discussion above for $\T^2$, in particular
Proposition~\ref{PROP:Gibbs1}, 
also holds on $\M$.
See~Section 4 of~\cite{OT}.
While the presentation in~\cite{OT} is given in the complex-valued setting, 
a straightforward modification yields the corresponding result 
for  the real-valued setting.

In the next subsection, we discuss the dynamical problem.
Our main goal in this paper is to construct dynamics
for the renormalized equation associated
with the Wick ordered Hamiltonian $H_\text{Wick}$
with initial data distributed according to the Gibbs measure $\Pk$.

\subsection{Dynamical problem: Wick ordered NLW}

We now  consider the following dynamical problem on $\T^2$
associated with the Wick ordered Hamiltonian: 
\begin{equation}
\begin{cases}
\displaystyle  \dt 
 \begin{pmatrix}
 u \\ v 
 \end{pmatrix}
 =  
 \begin{pmatrix}
 0& 1 \\ -1 & 0
 \end{pmatrix}
\frac{\dd H_{\text{Wick}}}{\dd(u, v )}\\
\rule[0mm]{0mm}{5mm}
(u, v) |_{t = 0} = (\phi_0^\o, \phi_1^\o),
\end{cases}
\label{Hamil3}
\end{equation}

\noi
where the initial data $(\phi_0^\o, \phi_1^\o)$ is distributed according to the Gibbs measure $\Pk$.
In view of the absolute continuity of $\Pk$ with respect to the Gaussian measure $\mu$
(Proposition~\ref{PROP:Gibbs1}), 
we consider 
the random initial data $(\phi_0^\o, \phi_1^\o)$
distributed according to $\mu$ in the following discussion.
Namely, we assume that 
\begin{align}
 \big(\phi_0^\o,  \phi_1^\o\big)
=\bigg( \sum_{n \in \Z^2} \frac{g_{0, n}(\o)}{\jb{n}_\rho}e^{in\cdot x},
\sum_{n \in \Z^2} g_{1, n}(\o)e^{in\cdot x}\bigg),
\label{Gauss1}
\end{align}

\noi
where $\{g_{0, n}, g_{1, n}\}_{n \in \Z^2}$
is as in \eqref{G5}.
Note that, at this point, the  potential part 
$ \frac1{2m+2} \int_{\T^2} :\! u^{2m+2}\!: dx$
of the Wick ordered Hamiltonian 
is defined  only for $u$ distributed according to the Gaussian measure $\mu$
via \eqref{Wick2}.
In the following, we extend this definition to a wider class of functions
in order to treat the Cauchy problem~\eqref{Hamil3}.

Given $N \in \N$, 
define 
the truncated Wick  ordered Hamiltonian $H^N_{\text{Wick}}$ by 
\begin{align}
H^N_\text{\tiny Wick}(u, v) =\frac{1}{2}\int_{\T^2}\,\big(\rho u^2 +  |\nb u|^2\big)\, dx
+ \frac{1}{2}\int_{\T^2} v^2dx
+ \frac1{2m+2} \int_{\T^2} :\! (\P_N u)^{2m+2}\!: dx
\label{Hamil4}
\end{align}

\noi
and consider the associated Hamiltonian dynamics:
\begin{equation*}
\begin{cases}
\displaystyle  \dt 
 \begin{pmatrix}
 u_N \\ v_N 
 \end{pmatrix}
 =  
 \begin{pmatrix}
 0& 1 \\ -1 & 0
 \end{pmatrix}
\frac{\dd H^N_{\text{Wick}}}{\dd(u_N, v_N )}\\
\rule[0mm]{0mm}{5mm}
(u_N, v_N) |_{t = 0} = (\phi_0^\o, \phi_1^\o).
\end{cases}
\end{equation*}

\noi
Thanks to 
 \eqref{Wick1}
and $\dx H_k(x; \s) = k H_{k-1}(x; \s)$, 
we can rewrite the system \eqref{Hamil4} as the
following truncated  Wick ordered NLW: 
\begin{align}
\begin{cases}
\dt^2 u_N -  \Dl  u_N + \rho u_N  +  \P_N\big[:\! (\P_N u_N)^{2m+1} \!:\big] = 0 \\
(u_N, \dt u_N) |_{t = 0} = ( \phi_0^\o,  \phi_1^\o),
\end{cases}
\label{WNLW1}
\end{align}

\noi
where the truncated Wick ordered nonlinearity is interpreted as
\begin{align*}
\P_N\big[:\! (\P_N u_N)^{2m+1} \!:\big]
= \P_N \big[ H_{2m+1}(\P_N u_N; \s_N)\big].
\end{align*}

Let $ z = z^\o$ denote the random linear solution:
\begin{align}
z(t) = S(t) (\phi_0^\o,  \phi_1^\o)
= 
\cos (t \jb{\nb}_\rho) \phi_0^\o
+ \frac{\sin (t \jb{\nb}_\rho)}{\jb{\nb}_\rho}   \phi_1^\o,
\label{lin1}
\end{align}

\noi
where $\jb{\nb}_\rho = \sqrt{\rho - \Dl}$.
In view of the Duhamel formula, it is natural to decompose the solution $ u_N$ to \eqref{WNLW1}
as 
\begin{align*}
 u_N = z + w_N.
\end{align*}

\noi
Note that we have $\P_N w_N = w_N$.
By  recalling 
the following  identities for the Hermite polynomials:
\begin{align}
H_k(x+y) & = \sum_{\l = 0}^k
\begin{pmatrix}
k \\ \l
\end{pmatrix}
H_\l(y)\cdot x^{k - \l} 
\qquad \text{and}\qquad 
H_k(x; \s )  = \s^\frac{k}{2} H_k(\s^{-\frac{1}{2}} x), 
\label{Herm1}
\end{align}

\noi
we have 
\begin{align}
:\! (\P_N u_N)^{2m+1}\!\!:  \, 
& = H_{2m+1}( z_N +  w_N; \s_N) \notag\\
&   = \sum_{\l = 0}^{2m+1}
\begin{pmatrix}
2m +1 \\ \l
\end{pmatrix}
 H_\l(z_N; \s_N)\cdot   w_N^{2m+1 - \l}, 
\label{Herm2}
\end{align}

\noi
where
$ z_N = \P_N z$.
This shows that  applying the Wick ordering to the monomial 
\begin{equation}
(\P_N u_N)^{2m+1} = 
 (z_N +  w_N)^{2m+1}
= \sum_{\l=0}^{2m+1} {2m+1\choose \ell} z_N^\l\cdot   w_N^{2m+1-\l}
\label{Herm3}
\end{equation}

\noi
is equivalent to 
Wick ordering all the monomials $z_N^\l$.
Namely, replacing each $z^\l$ in \eqref{Herm3} by 
\[:\! z_N^\l\!: \, =  H_\l(z_N; \s_N)\]

\noi
yields the Wick ordered monomial $:\! (\P_N u_N)^{2m+1}\!\!:  $ via 
 \eqref{Herm2}.
In Proposition~\ref{PROP:PStr} below, 
we prove that 
\[:\! z_N^\l\!:\, \in L^p(\O; L^q([-T, T]; W^{-\eps, r}(\T^2)) )\]

\noi
for any $p, q, r  < \infty$, $T> 0$, and $\eps > 0$
with a  bound uniform in $N$.
Moreover, the sequence $\big\{:\! z_N^\l\!:\big\}_{N \in \N}$ is a Cauchy sequence  in the same space, 
thus allowing us to define 
\begin{align}
:\! z^\l\!:\
= \ :\! z_\infty^\l\!:\, \stackrel{\text{def}}{=} \lim_{N \to \infty}
:\! z_N^\l\!:
\label{Wick4}
\end{align}

\noi
in $ L^p(\O; L^q([-T, T]; W^{-\eps, r}(\T^2)) )$
for any $p, q, r  < \infty$, $T > 0$,  and $\eps > 0$
(and for any $\l \in \N$).
Now, consider a function $u$ of the form
\begin{align}
u = z + w
\label{decomp2}
\end{align}

\noi
 for some ``nice'' $w$.
Then, 
we can use  \eqref{Herm1} and \eqref{Wick4}
to  define  the Wick ordered monomial 
$:\! u^{2m+1}\!:$ for functions $u$ of the form \eqref{decomp2} by 
\begin{align}
:\!  u^{2m+1}\!\!:  \, 
& 
= \ :\!  (z + w)^{2m+1}\!\!:  \, 
  = \sum_{\l = 0}^{2m+1}
\begin{pmatrix}
2m +1 \\ \l
\end{pmatrix}
:\! z^\l\!: \cdot  \, w^{2m+1 - \l}.
\label{Herm4}
\end{align}

\noi
Hence, we finally arrive at the 
{\it defocusing Wick ordered NLW}: 
\begin{align}
\begin{cases}
\dt^2 u -  \Dl  u + \rho u \, +  :\! u^{2m+1} \!:\, = 0 \\
(u, \dt u) |_{t = 0} = (\phi_0^\o, \phi_1^\o), 
\end{cases}
\label{WNLW2}
\end{align}

\noi
where $ (\phi_0^\o, \phi_1^\o)$ is as in \eqref{Gauss1}.

Before we state our main result, 
we first recall two critical regularities associated
with \eqref{NLW1} on $\R^2$ with $\rho = 0$.
On the one hand,  the scaling symmetry for \eqref{NLW1}
 induces the so-called scaling critical Sobolev index:
$s_1 = 1 - \frac 1m$.
On the other hand, 
the Lorentzian invariance (conformal symmetry)
induces another critical regularity: 
$s_2 = \frac{3}{4} - \frac{1}{2m}$
(at least in the focusing case).
Hence, we set $s_\text{crit}$ by 
\begin{align*}
s_\text{crit} = \max\bigg(1 - \frac{1}{m}, \frac 34 - \frac 1{2m}\bigg)= 
\begin{cases}
\frac14 \quad &\text{if} \quad  m=1, \\
1-\frac1m \;&  \text{if} \quad m\geq 2. 
\end{cases}
\end{align*}

\noi
We now state our main result.

\begin{theorem}\label{THM:LWP}
Let  $\M = \T^2$, $m \in \N$, and $\rho > 0$.
Then, the  Wick ordered NLW \eqref{WNLW2}
is almost surely locally well-posed
with respect to 
the Gaussian measure $\mu$ defined in \eqref{G3}.
More precisely,
letting 
 $(\phi_0^\o, \phi_1^\o)$ be as in \eqref{Gauss1}, 
there exist $ C,  c>0$ such that
for each $T\ll 1$,
there exists a set $\O_T  \subset \O$ with the following properties:

\smallskip
\begin{itemize}
\item[\textup{(i)}]
$P(\O_T^c) \leq C \exp\big(-\frac{1}{T^{c} }\big)$,

\smallskip
\item[\textup{(ii)}]
For each $\o \in \O_T$, there exists a (unique) solution $u$
to \eqref{WNLW2}
with $(u, \dt u) |_{t = 0} = (\phi_0^\o, \phi_1^\o)$ 
in the class
\[ S(t)( \phi_0^\o,  \phi_1^\o) + C([-T, T]; H^{s} (\T^2)) 
\cap X^{s, \frac 12 +}_T
\subset C([-T, T];H^{-\eps}(\T^2))\]

\noi
for any $s \in ( s_\textup{crit}, 1)$ and $\eps > 0$.
Here, $X^{s, \frac 12 +}_T$ denotes the local-in-time version 
of the hyperbolic Sobolev space.
See Section \ref{SEC:LWP}.
\end{itemize}

\end{theorem}

\noi

We emphasize that 
the Wick ordered NLW \eqref{WNLW2} is defined {\it only} for functions 
$u$ of the form~\eqref{decomp2}.
Then, the residual term $w = u - z$ satisfies
the following perturbed Wick ordered NLW:
\begin{align}
\begin{cases}
\dt^2 w -  \Dl  w + \rho w  \, +  :\! (w+z) ^{2m+1} \!:\ = 0 \\
(w, \dt w) |_{t = 0} = (0, 0).
\end{cases}
\label{WNLW3}
\end{align}

\noi
By writing \eqref{WNLW3} 
in the Duhamel formulation, 
we obtain 
\begin{align}
w(t)
& = 
- \int_{0}^t  \frac{\sin ((t-t')\jb{\nb}_\rho)}{\jb{\nb}_\rho}  
 :\! (w+z) ^{2m+1} (t')\!: dt'\notag\\
& = 
-  \sum_{\l = 0}^{2m+1}
\int_{0}^t  \frac{\sin ((t-t')\jb{\nb}_\rho)}{\jb{\nb}_\rho}  
\begin{pmatrix}
2m +1 \\ \l
\end{pmatrix}
:\! z^\l(t')\!: \cdot  \, w^{2m+1 - \l}(t')dt'.
\label{WNLW4}
\end{align}

\noi
We prove Theorem~\ref{THM:LWP}
by solving the  fixed point problem \eqref{WNLW4}
for $w$ 
in $C([-T, T]; H^{s} (\T^2))\cap X^{s, \frac 12 +}_T$,  $s > s_\textup{crit}$.
In Section \ref{SEC:sto}, we study the regularity of the random linear solution $z$
and the associated Wick ordered monomials $:\! z^\l \!:$.
In particular, while they are rough, 
$:\! z^\l \!:$ enjoys enhanced integrability both in space and time.
See Proposition~\ref{PROP:PStr}.
In Section \ref{SEC:LWP}, 
we then use the standard Fourier restriction norm method
to solve the fixed point problem \eqref{WNLW4}.
The original idea of this argument with the decomposition \eqref{decomp2}
 appears in McKean \cite{McKean} and 
Bourgain \cite{BO96} in the context of 
the nonlinear Schr\"odinger equations on $\T^d$, $d = 1, 2$.
 See also Burq-Tzvetkov \cite{BT1}.
In the field of the  stochastic PDEs, this method 
is  known  as 
 Da Prato-Debussche trick \cite{DPD}.

\begin{remark}\rm
As in the study of singular stochastic PDEs, 
our proof consists of factorizing the ill-defined solution map:
$(\phi_0^\o, \phi_1^\o)\mapsto u$
into a canonical lift followed by a (continuous) solutions map $\Psi$: 
\begin{align*}
(\phi_0^\o, \phi_1^\o)
\stackrel{\text{lift}}{\longmapsto} (z_1^\o, z_3^\o, \dots, z_{2m+1}^\o ) \stackrel{\Psi}{\longmapsto}
&  \, w \in C([-T, T]; H^s(\T^2))\notag\\
\longmapsto &  \, u =  z + w \in C([-T, T]; H^{-\eps}(\T^2)), 
\end{align*}

\noi 
for $s \in ( s_\text{crit}, 1)$ and $\eps > 0$, 
where $z_k \stackrel{\text{def}}{=} \, :\!z^k\!\!:\, $.
On the one hand, we use probability theory to construct the data set $\{z_{2j+1}\}_{j = 0}^m$
in the first step.
On the other hand, the second step is entirely deterministic.
Moreover, the solution map $\Psi$ in the second step is continuous
from $\prod_{j = 0}^m S^j_T$
to $X^{s, \frac{1}{2}+}_T$, 
where $S^j_T$ denotes some appropriate Strichartz space
for $z_{2j+1}$.  See Section \ref{SEC:LWP}.

\end{remark}

\begin{remark}\label{REM:approx}\rm
The same almost sure local well-posedness
holds for the truncated Wick ordered NLW \eqref{WNLW1}.
More precisely,  we can choose $\O_T$, independent of $N \in \N$, such that the statement in Theorem~\ref{THM:LWP}
holds for \eqref{WNLW2} and \eqref{WNLW1}.
Moreover, by possibly shrinking the time, one can also prove that 
the solution $u_N = u_N^\o$ to \eqref{WNLW1}
converges to 
the solution $u = u^\o$ to \eqref{WNLW2} as $N \to \infty$.

\end{remark}

Once we have almost sure local well-posedness of \eqref{WNLW2}, 
the invariant measure argument by 
Bourgain \cite{BO94, BO96}
yields the following almost sure global well-posedness
of \eqref{WNLW2}
and invariance of the Gibbs measure $\Pk$.

\begin{theorem}\label{THM:GWP}

Let  $\M = \T^2$, $m \in \N$, and $\rho > 0$.
Then, the  defocusing Wick ordered NLW \eqref{WNLW2}
is almost surely globally  well-posed
with respect to 
the Gibbs measure $\Pk$.
Moreover, $\Pk$ is invariant under the dynamics of 
\eqref{WNLW2}.

\end{theorem}

The proof of Theorem~\ref{THM:GWP}
exploits the invariance of the truncated Gibbs measure $\Pkn$
for the  truncated  Wick ordered NLW \eqref{WNLW1}
and combines it with an approximation argument.
See Remark \ref{REM:approx}.
As this argument is standard by now, we omit the proof.
See Bourgain \cite{BO96} and  Burq-Tzvetkov \cite{BT2}
for details.

\begin{remark}\rm
We point that 
the convergence result in Remark \ref{REM:approx}
and invariance of the Gibbs measure in Theorem \ref{THM:GWP}
already appear (without a proof) in the lecture note by Bourgain \cite{BoPCMI}.
See \cite[Theorem 111 on p.\,63]{BoPCMI}
and a comment that follows (118) on p.\,64 in \cite{BoPCMI}.
To the best of our knowledge, however, 
there seems to be no proof available
in a published paper.
In fact, one of the main purposes
of this paper is to present the details of the proof of Bourgain's claim in \cite{BoPCMI}.
\end{remark}

\medskip

Next, we briefly discuss the situation 
when the spatial domain $\M$
is  a two-dimensional compact Riemannian manifold without boundary or a bounded domain in $\R^2$
(with the Dirichlet or Neumann boundary condition). 
In this case, one can 
 exploit the invariance
of the truncated Gibbs measures $\Pkn$
for~\eqref{WNLW1}
 to construct global-in-time {\it weak} solutions (without uniqueness) 
 to the Wick ordered NLW~\eqref{WNLW2}.
Moreover, it also allows us to establish 
 invariance of the Gibbs measure $\Pk$ in some mild sense.

\begin{theorem}\label{THM:GWP2}
Let $m \in \N$ and $\rho > 0$.  
Let $\M$
be  a two-dimensional compact Riemannian manifold without boundary or a bounded domain in $\R^2$
(with the Dirichlet or Neumann boundary condition). 
In the latter case with the Dirichlet boundary condition, 
we can also take $\rho = 0$.
Then, there exists a set $\Si$ of full measure
with respect to $\Pk$
such that for every $\phi \in \Si$, 
the  defocusing  Wick ordered NLW~\eqref{WNLW2} with 
initial data distributed according to $\Pk$
has a global-in-time solution 
$ u \in C(\R; H^s(\M))$
for any $s < 0$.
Moreover, for all $t \in \R$, 
the law of the random function $(u, \dt u)(t)$ is given by~$\Pk$.
\end{theorem}

In \cite{OT}, we proved an analogous result 
for the defocusing Wick ordered nonlinear Schr\"odinger equations
on $\M$.   Theorem~\ref{THM:GWP2} follows from  repeating 
the argument presented in~\cite{OT} with systematic modifications
and thus we omit details.
See also  \cite{AC, DPD, BTT1, ORT}.
The main  ingredient for Theorem~\ref{THM:GWP2}
is to establish tightness (= compactness)
of measures $\nu_N$ on space-time functions, emanating
from the truncated Gibbs measure $\Pkn $
and then upgrading the weak convergence of $\nu_N$ (up to a subsequence)
to an almost sure convergence of the corresponding random variables 
via  Skorokhod's theorem.
Due to the compactness argument, 
Theorem~\ref{THM:GWP2} claims 
 only the existence of a global-in-time solution $u$. 
Lastly, note that
Theorem~\ref{THM:GWP2} only claims
that the law  of the $\H^s$-valued random variable $(u, \dt u)(t)$
is given by the Gibbs measure $\Pk$
for any $t \in \R$.
In particular,  
this mild invariance for a general geometric setting is weaker than the invariance stated in Theorem~\ref{THM:GWP} for the Wick ordered NLW \eqref{WNLW2} on $\T^2$.

\begin{remark}\rm
On the one hand, 
the defocusing/focusing nature of the equation does not play any role
in the 
almost sure local well-posedness result (Theorem~\ref{THM:LWP})
and thus Theorem~\ref{THM:LWP} also holds
 in the focusing case.  It can also be extended to Wick ordered even power monomials
 in the equation.
On the other hand, 
the defocusing nature of the equation 
plays a crucial role in the proof of Proposition~\ref{PROP:Gibbs1}
and hence in Theorems~\ref{THM:GWP} and~\ref{THM:GWP2}.
In the focusing case (i.e.~with $- u^{2m+1}$, $m \in \N$, in \eqref{NLW1}), it is known that the Gibbs measure 
can not be normalized in the two dimensional case.
See Brydges-Slade \cite{BS}.
Lastly, we point out
that in the case of the quadratic nonlinearity (which is neither defocusing nor focusing), 
one can introduce 
the following modified Gibbs measure:
\[ d P_2^{(3)} = Z^{-1} e^{-\frac 13 \int : u^3: \,  - A (\int : u^2:\,)^2} d\mu\]

\noi
for sufficiently large $A \gg 1$
and study the associated dynamical problem.
See \cite{BoPCMI} for the construction of this modified Gibbs measure $P_2^{(3)}$.

\end{remark}

\subsection{Wick ordered NLW as a scaling limit}

As an application of the local well-posedness argument, 
we  show how the Wick ordered NLW \eqref{WNLW2}
appears as a scaling limit
of non-renormalized NLW equations on dilated tori.
This part of the discussion is strongly motivated
by the weak universality
result for the  Wick ordered stochastic NLW
on $\T^2$ studied  by the first author with Gubinelli and Koch
in \cite{GKO}.

Fix $\rho > 0$.
 Given small $\eps > 0$, we consider the following  non-renormalized NLW equation 
on a dilated torus $\T_\eps^2 \stackrel{\text{def}}= (\eps^{- 1} \T)^2$:
\begin{align}
\begin{cases}
 \dt^2 v_\eps - \Dl v_\eps + \rho_\eps v_\eps  = f (v_\eps) \\
 (v_\eps, \dt v_\eps) |_{t = 0}=  (\psi_{\eps, 0}^\o, \psi_{\eps, 1}^\o),
\end{cases}
\qquad 
(t,x)\in  \R \times \T_\eps^2 
\label{XNLW1}
\end{align}

\noi
with Gaussian random initial data $(\psi_{\eps, 0}^\o, \psi_{\eps, 1}^\o)$, 
where 
 $f : \mathbb{R} \rightarrow
\mathbb{R}$ is a  smooth odd\footnote{It follows from the proof of Theorem~\ref{THM:WNLW} 
that it suffices to assume that  $f(0) = f''(0) = 0$
 for the cubic case considered in Theorem~\ref{THM:WNLW}.}
function with the following bound: 
\begin{align}
|f^{(4)}(x)| \les 1 + |x|^{M}
\label{bound}
\end{align}

\noi
for some $M \geq 0$
and $\rho_\eps$ is a parameter to be chosen later.
In the following, we choose 
$\psi_{\eps, 0}^\o$ and 
$\psi_{\eps, 1}^\o$ to be
a smoothed Ornstein-Uhlenbeck process
and  a smoothed white noise on $\T_\eps^2$, respectively. 
For the sake of concreteness, we set\footnote{Note that 
$\big\{\eps e^{i n \cdot x} \big\}_{n \in (\eps \Z)^2}$ forms an orthonormal basis of $L^2(\T_\eps^2)$.
Moreover, recall that the Fourier-Wiener series 
\[\sum_{n \in (\eps \Z)^2}
\frac{  g_{0, n}}{\eps^{-1}|n|} \frac{e^{i n \cdot x}}{\eps^{-1}}\] 

\noi
represents the periodic Wiener process on $\T_\eps^2$.
} 
\begin{align*}
(\psi_{\eps, 0}^\o, 
\psi_{\eps, 1}^\o )= 
 \Bigg( \sum_{\substack{n \in (\eps \Z)^2\\|n|\leq 1}} 
\frac{  g_{0, \eps^{-1} n}}{\eps^{-1}\sqrt{\eps^2 \rho + |n|^2}} \frac{e^{i n \cdot x}}{\eps^{-1}}, 
 \sum_{\substack{n \in (\eps \Z)^2\\|n|\leq 1}}    g_{1, \eps^{-1} n} \frac{e^{i n \cdot x}}{\eps^{-1}}\Bigg), 
\end{align*}

\noi
where $\{g_{0, n}, g_{1, n}\}_{n \in \Z^2}$ is as in \eqref{G5}.
Our main goal is to study  the behavior of the solution to  \eqref{XNLW1} as $\eps \to 0$
by applying a suitable scaling.

Let $u_{\eps} (t,x) \overset{\text{def}}{=} \eps^{- 1} v_\eps (\eps^{-1} x, \eps^{-1} t)$.
Then,   $u_{\eps}$ satisfies
\begin{align}
\begin{cases}
 \dt^2 u_\eps - \Dl u_\eps + \rho u_\eps
  = \eps^{ - 3} \big\{f
(\eps u_{\eps}) + \eps (\eps^2 \rho -  \rho_\eps) u_{\eps}\big\}\\
(u_\eps, \dt u_\eps) |_{t = 0} = (\phi^\o_{\eps, 0}, \phi^\o_{\eps, 1}),
\end{cases}
\qquad 
(t,x)\in  \R \times \T^2,  
\label{XNLW2}
\end{align}

\noi
where $(\phi^\o_{\eps, 0}, \phi^\o_{\eps, 1})$ on $\T^2$ is given by 
\begin{align}
(\phi_{\eps, 0}^\o, 
\phi_{\eps, 1}^\o )= 
  \Bigg( \sum_{\substack{n \in \Z^2\\|n|\leq \eps^{-1}}} 
\frac{  g_{0, n}}{\sqrt{ \rho + |n|^2}} e^{i n \cdot x}, 
 \sum_{\substack{n \in \Z^2\\|n|\leq \eps^{-1}}}    g_{1, n}e^{i n \cdot x}\Bigg).
\label{A2}
\end{align}

\noi
Note that 
$(\phi^\o_{\eps, 0}, \phi^\o_{\eps, 1})$ converges  to 
$(\phi^\o_{0}, \phi^\o_{1})$ in \eqref{Gauss1}
distributed according to the Gaussian measure $\mu$ defined in \eqref{G3}.

The rescaled model \eqref{XNLW2} on $\T^2$ allows us to study 
 the large temporal and spatial  scale behavior of
the solution $v_\eps$ to \eqref{XNLW1}.
Moreover, by suitably choosing $\rho_\eps$, 
the family $\{u_\eps\}_{\eps>0}$ converges
to the solution $u$ to the Wick ordered NLW on $\T^2$
with a parameter $\ld = \ld(f)$,
depending only on $f$.

\begin{theorem}\label{THM:WNLW}
Let $\rho > 0$.  Then,  there exists a choice of $\rho_\eps$ such that, as $\eps \to 0$, 
 the family of 
the solutions   $\{u_{\eps}\}_{\eps>0}$ to \eqref{XNLW2} converges 
almost surely 
 to  the solution $u$ to the following  Wick ordered cubic NLW on $\T^2$:
  \begin{equation}
\begin{cases}
\dt^2 u - \Dl  u + \rho u = \lambda : \!u^3 \!: \\
(u, \dt u) |_{t = 0} = (\phi^\o_{0}, \phi^\o_{1}),
\end{cases}
    \label{NLW_ld}
  \end{equation}
where the convergence takes place  
in $C([-T_\o, T_\o]; H^s(\T^2))$, $s < 0$,  for some  $T_\o>0$.
   Here, the constant $\lambda$ is given by  $\ld =\frac{f^{(3)} (0)}{6} $, 
   depending only on the function $f$.

\end{theorem}

This theorem shows  a kind of weak universality  for the Wick ordered NLW.
See \cite{GKO} for a similar result for the Wick ordered stochastic NLW.
We also refer readers to 
\cite{HQ, GP1, GP2} for more discussion on weak universality 
(for stochastic parabolic equations, in particular the KPZ equation).

\begin{remark}\rm
By starting with   the following  NLW on $\T_\eps^2$: 
\begin{align*}
\begin{cases}
 \dt^2 v_\eps - \Dl v_\eps + \rho_\eps v_\eps  + \sum_{j = 1}^{m-1}
a_j (\eps)v_\eps^{2j+1}
 = f (v_\eps) \\
 (v_\eps, \dt v_\eps) |_{t = 0}=  (\psi_{\eps, 0}^\o, \psi_{\eps, 1}^\o),
\end{cases}
\end{align*}

\noi
  we can tune the $m$ parameters $\rho_\eps, a_j(\eps)$, $j = 1, \dots, m-1$, 
  such that by a small modification of the proof of Theorem~\ref{THM:WNLW},  
we obtain the following Wick ordered NLW:
    \begin{equation*}
    \begin{cases}
\dt^2 u - \Dl  u + \rho u= \lambda : \!u^{2m+1} \!: \\
(u, \dt u) |_{t = 0} = (\phi^\o_{0}, \phi^\o_{1}),
\end{cases}
  \end{equation*}

\noi
for some $\ld = \ld(f)$, as $\eps \to 0$.
  In this case, 
  one needs to use the scaling 
   $u_{\eps} (t,x) = \eps^{- \g} v_\eps (\eps^{-1} x, \eps^{-1} t)$
   for some suitably chosen $\g = \g(m)>0$
   and also assume a bound analogous to \eqref{bound}
   for a higher order derivative of $f$.

\end{remark}

\section{Probabilistic tools}
\label{SEC:sto}

In this section, we first recall basic probabilistic tools.
Then, we prove a uniform (in $N$) bound on the Wick ordered monomials 
$:\! z_N^\l\!:\, = H_\l(z_N, \s_N)$, consisting of 
the random linear solution (Proposition~\ref{PROP:PStr}).
Moreover, we prove that $\big\{:\! z_N^\l\!:\big\}_{N \in \N}$ is a Cauchy sequence,
allowing us to define
$:\! z^\l\!:$ by \eqref{Wick4}.

\subsection{Hermite polynomials
and white noise functional}
First, recall the Hermite polynomials $H_k(x; \s)$
defined via the generating function:
\begin{equation}
F(t, x; \s) : =  e^{tx - \frac{1}{2}\s t^2} = \sum_{k = 0}^\infty \frac{t^k}{k!} H_k(x;\s).
\label{H1}
 \end{equation}
	
\noi
For simplicity, we set 
$H_k(x) : = H_k(x; 1)$.
In the following, we list the first few Hermite polynomials
for readers' convenience:
\begin{align}
\begin{split}
& H_0(x; \s) = 1, 
\qquad 
H_1(x; \s) = x, 
\qquad
H_2(x; \s) = x^2 - \s, \\
& H_3(x; \s) = x^3 - 3\s x, 
\qquad 
H_4(x; \s) = x^4 - 6\s x^2 +3\s^2.
\end{split}
\label{H1a}
\end{align}

Next, we define the white noise functional.
Let  $\xi(x;\o)$ be the (real-valued) mean-zero Gaussian white noise on $\T^2$
defined by
\[ \xi(x;\o) = \sum_{n\in \Z^2} g_n(\o) e^{in\cdot x},\]

\noi
where 
$\{g_{ n} \}_{n \in \Z^2}$
is a sequence of independent standard complex-valued Gaussian
random variables
conditioned that $g_{ -n} = \cj{g_{n}}$, $n \in \Z^2$.
It is easy to see that $\xi \in \H^s(\T^2) \setminus \H^{-1}(\T^2)$, $s < -1$, almost surely.
In particular, $\xi$ is a distribution, acting  on smooth functions.
In fact, the action of $\xi$ can be defined on $L^2(\T^2)$.

 We define the white noise functional  $W_{(\cdot)}: L^2(\T^2) \to L^2(\O)$
 by 
\begin{equation}
 W_f (\o) = \jb{f, \xi(\o)}_{L^2} = \sum_{n \in \Z^2} \ft f(n) \cj{g_n}(\o)
\label{W0}
 \end{equation}

\noi
for a real-valued function $f \in L^2(\T^2)$.
Note that $W_f=\xi(f)$ is basically the Wiener integral of $f$.
In particular, 
$W_f$ is a real-valued Gaussian random variable
with mean 0 and variance $\|f\|_{L^2}^2$.
Moreover, 
$W_{(\cdot)}$ is unitary:
\begin{align}
 E\big[ W_f W_h ] = \jb{f, h}_{L^2_x}
 \label{W0a}
\end{align}

\noi
for $f, h \in L^2(\T^2)$.
The following lemma extends the relation \eqref{W0a}
to a more general setting.

\begin{lemma}\label{LEM:W1}
Let $f, h \in L^2(\T^2)$ such that $\|f\|_{L^2} = \|h\|_{L^2} = 1$.
Then, for $k, m \in \Z_{\geq 0}$, we have 
\begin{align*}
\E\big[ H_k(W_f)H_m(W_h)\big]
=  \dl_{km} k! [\jb{f, h}]^k.
\end{align*}

\noi
Here, $\dl_{km}$ denotes the Kronecker's delta function.
\end{lemma}

\noi
This lemma  follows from 
computing the left-hand side of 
\[\E [ F(t, W_f) F(s, W_h)]
=  \sum_{k, m  = 0}^\infty \frac{t^k}{k!}
 \frac{s^m}{m!}
\E\big[ H_k(W_f)H_m(W_h)\big]\]

\noi
and comparing the coefficients.
See \cite{DPT1, ORT} for details.

We also recall the following Wiener chaos estimate
\cite[Theorem I.22]{Simon}.

\begin{lemma}\label{LEM:hyp3}
Fix $k \in \mathbb{N}$ and $c(n_1, \dots, n_k) \in \C$.
Given 	
 $d \in \mathbb{N}$, 
 let $\{ g_n\}_{n = 1}^d$ be 
 a sequence of  independent standard complex-valued Gaussian random variables
 and set $g_{-n} = \cj{g_n}$.
Define $S_k(\o)$ by 
\begin{align*}
 S_k(\o) = \sum_{\G(k, d)} c(n_1, \dots, n_k) g_{n_1} (\o)\cdots g_{n_k}(\o),
 \end{align*}

\noi
where $\G(k, d)$ is defined by
\[ \G(k, d) = \big\{ (n_1, \dots, n_k) \in \{0, \pm1, \dots,\pm d\}^k \big\}.\]

\noi
Then, for $p \geq 2$, we have
\begin{equation}
 \|S_k \|_{L^p(\O)} \leq (p-1)^\frac{k}{2}\|S_k\|_{L^2(\O)}.
\label{hyp4}
 \end{equation}

\end{lemma}

The crucial point is that  the constant in \eqref{hyp4} 
is {\it independent} of $d \in \N$.
This lemma is a direct corollary to the
  hypercontractivity of the Ornstein-Uhlenbeck
semigroup due to Nelson \cite{Nelson2}.

\subsection{Stochastic estimate on Wick ordered monomials}

In this subsection, we study 
the Wick ordered monomials
$:\! z_N^\l\!:$ and 
$:\! z^\l\!:$,
 consisting of the random linear solution $z$ defined in \eqref{lin1}.
From \eqref{Gauss1} and \eqref{lin1}, we have
\begin{align}
\ft z(t, n) = 
\frac{\cos (t \jb{n}_\rho)}{\jb{n}_\rho} g_{0, n}
+ \frac{\sin (t \jb{n}_\rho)}{\jb{n}_\rho}   g_{1, n}.
\label{PStr1}
\end{align}

\noi
In order to avoid the combinatorial complexity in higher ordered monomials, 
we use  the white noise functional as in \cite{OT}.
We, however, need to adapt the white noise functional to $z(t)$.
In view of \eqref{PStr1}, 
we define  the white noise functional 
$W^t_{(\cdot)}: L^2(\T^2) \to L^2(\O)$
with a parameter $t \in \R$ by 
\begin{equation}
 W^t_f (\o) = \jb{f, \xi^t(\o)}_{L^2} = \sum_{n \in \Z^2} \ft f(n) \cj{g^t_n}(\o).
\label{PStr2}
 \end{equation}

\noi
Here,  $\xi^t$ denotes (a specific realization of) the white noise on $\T^2$ given by 
\[\xi^t (x; \o) = \sum_{n\in \Z^2} g_n^t (\o) e^{in \cdot x},\]

\noi
where $g_n^t$ is define by 
\[g_n^t = \cos (t\jb{n}_\rho) g_{0, n}+\sin (t\jb{n}_\rho) g_{1, n}.\]

\noi
Note that $\E [g_n^t] = 0$ and 
$\text{Var}(g_n^t) =  
\cos^2 (t\jb{n}_\rho) +\sin^2 (t\jb{n}_\rho)  = 1$.
Thus,  for each fixed $t \in \R$, 
 $\{g_n^t\}_{n\in \Z^2}$ is a sequence of independent standard Gaussian random variables
 conditioned that $g^t_{-n} = \cj{g^t_n}$ for all $n \in \N$.
Therefore, the white noise functional
$W^t_{(\cdot)}$ defined in \eqref{PStr2} satisfies
the same properties as the standard white noise functional 
$W_{(\cdot)}$ defined in \eqref{W0}. 
Lastly, note that, in view of \eqref{PStr1},  
the random linear solution $z_N = \P_N z$ can be expressed as 
\begin{align}
z_N(t,x) =  \sum_{|n|\leq N} \frac{g_{n}^t(\o)}{\jb{n}_\rho}e^{in\cdot x}.
\label{PStr2a}
\end{align}

\noi
In the following, 
 we use the short-hand notation  $L_T^q= L^q_t([-T, T])$, etc.

\begin{proposition}\label{PROP:PStr}
Let  $\l \in \N$ and $\rho > 0$.
Then,
given $2\leq q, r <\infty$ and $\eps > 0$, 
there exist $C, c>0$ such that
\begin{align}
P\Big( \| \jb{\nb}^{-\eps} 
\! :\! z_N^\l\!:
\|_{L^q_T L^r_x}> \ld\Big)
\leq C\exp\bigg(-c\frac{\ld^\frac{2}{\l}}{T^{\frac{2}{q\l}}}\bigg)
\label{PStr3}
\end{align}
	
\noi
for any $T > 0$, $\ld > 0$, and any $N \in \N$.
Moreover, 
the sequence $\big\{:\! z_N^\l\!:\big\}_{N \in \N}$ is a Cauchy sequence  
in $L^p(\O; L^q([-T, T]; W^{-\eps, r}(\T^2)) )$.
In particular, denoting the limit by 
$:\! z^\l\!:$, we have
$:\! z^\l\!: \, \in L^q([-T, T]; W^{-\eps, r}(\T^2)) $ almost surely,
satisfying the tail estimate \eqref{PStr3}.
\end{proposition}

Before proceeding to the proof of Proposition~\ref{PROP:PStr}, we introduce some notations.
Let $\s_N$ be as in \eqref{G6}.
For {\it fixed} $x \in \T^2$ and 
 $N \in \N$, 
 we also define
\begin{align}
\eta_N(x) (\cdot) & \stackrel{\text{def}}{=} \frac{1}{\s_N^\frac{1}{2}}
\sum_{|n| \leq N} \frac{\cj{e_n(x)}}{\jb{n}_\rho}e_n(\cdot)
\qquad \text{and}\qquad 
\g_N (\cdot)  \stackrel{\text{def}}{=} 
\sum_{ |n| \leq N} \frac{1}{\jb{n}_\rho^2}e_n(\cdot),
\label{W3}
\end{align}
	
\noi
where $e_n(y) = e^{in\cdot y}$.
Note that $\eta_N(x)(\cdot)$ is real-valued
with 
$ \| \eta_N(x)\|_{L^2(\T^2)} = 1$
for all  $x \in \T^2$ and all $N \in \N$.
Moreover, we have 
\begin{align}
\jb{\eta_M(x), \eta_N(y)}_{L^2}
= \frac{1}{\s_M^\frac{1}{2}\s_N^\frac{1}{2}} \g_N(y-x)
= \frac{1}{\s_M^\frac{1}{2}\s_N^\frac{1}{2}} \g_N(x-y), 
\label{W4}
\end{align}

\noi
for fixed $x, y\in \T^2$
and  $M\geq N \geq 1$.

\begin{proof}

From \eqref{PStr2a}  and \eqref{W3}, we have 
\begin{align}
z_N (t,x) = \s_N^\frac{1}{2} \frac{z_N(t,x)}{\s_N^\frac{1}{2}}
= \s_N^\frac{1}{2}W^t_{\cj{\eta_N(x)}} 
= \s_N^\frac{1}{2}W^t_{\eta_N(x)} .
\label{PStr4}
\end{align}

\noi
Then, from \eqref{Herm1} and \eqref{PStr4}, we have
\begin{align}
:\! z_N^\l(t,x)\!:
\, = H_\l(z_N(t,x); \s_N ) =  \s_N^\frac{\l}{2}H_\l\big(W^t_{\eta_N(x)} \big).
\label{PStr5}
\end{align}

Given $n \in \Z^2$, 
define  $\G_\l(n)$ by 
\begin{align*}
\G_\l (n) = \{ (n_1, \dots, n_\l) \in (\Z^2)^\l: \, n_1 + \cdots + n_\l = n\}.
\end{align*}

\noi
Then, for $(n_1, \dots, n_\l) \in \G_\l(n)$, we have $\max_j |n_j| \ges |n|$.
Thus, 
it follows from  Lemma \ref{LEM:W1} with \eqref{PStr5} and \eqref{W4} that
\begin{align}
\| \jb{ \, :\! z_N^\l(t)\!:, e_n} \|_{L^2(\O)}^2  
& = \s_N^\l \int_{\T^2_x \times \T^2_y}\cj{e_n (x)} e_n(y) 
 \int_\O H_\l\big(W^t_{\eta_N(x)} \big)
\cj{H_\l\big(W^t_{\eta_N(y)} \big)}
 dP dx dy \notag\\
& = \l! \int_{\T^2_x \times \T^2_y} [\g_N(x-y)]^\l  \cj{e_n( x-y)}dx dy \notag \\
& = \l! \cdot \F [\g_N^\l ](n) 
 = \l ! \sum_{\substack{\G_\l(n)\\ |n_j| \leq N}}
\prod_{j= 1}^\l \frac{1}{\jb{n_j}_\rho^2}
\les 
\frac{1}{\jb{n}^{2(1-\theta)}}
\label{PStr5a}
\end{align}

\noi
for any $\theta > 0$.
On the other hand, for $n \ne n'$, we have 
\begin{align}
\int_\O  \jb{ \, :\! z_N^\l(t)\!:, e_n} & \cj{  \jb{ \, :\! z_N^\l(t)\!:, e_{n'}}}  dP \notag\\
& = \s_N^\l \int_{\T^2_x \times \T^2_y}\cj{e_n(x) }e_{ n'} (y) 
 \int_\O H_\l\big(W^t_{\eta_N(x)} \big)
\cj{H_\l\big(W^t_{\eta_N(y)} \big)}
 dP dx dy \notag\\
& = \l! \int_{\T^2_x} \int_{\T^2_y}[\g_N(x-y)]^\l 
\cj{e_{ n'} (x-y)} dy \, 
\cj{e_n(x)} e_{n'}(x) dx \notag \\
& = \l! \cdot \F [\g_N^\l ](n') 
 \int_{\T^2_x}
\cj{e_n(x)} e_{n'}(x) dx 
= 0.
\label{PStr5b}
\end{align}

\noi
Hence, given  $x \in \T^2$ and $t \in \R$, it follows
from \eqref{PStr5a} and \eqref{PStr5b} that 
\begin{align}
\|\jb{\nb}^{-\eps}
 :\! z_N^\l(t,x)\!: \|_{L^2(\O)}
& =  \bigg\|\sum_{n \in \Z^2} \jb{n}^{-\eps} 
\F_x\big[ \, :\! z_N^\l(t)\!:\big](n) \, e^{in \cdot x} \bigg\|_{L^2(\O)} \notag \\
& \leq C_\l \bigg(\sum_{n \in \Z^2} \jb{n}^{-2\eps} 
\F [\g_N^\l ](n) \bigg)^\frac{1}{2}
\leq  
 C_\l \bigg(\sum_{n \in \Z^2} \jb{n}^{-2(1+\eps -\theta)} \bigg)^\frac{1}{2} \notag \\
&  < \infty, 
\label{PStr5c}
\end{align}

\noi
uniformly in $N \in \N$, 
as long as $0<  \theta < \eps$.

Fix  $2\leq q, r <\infty$.
Then, by Minkowski's integral inequality, Lemma \ref{LEM:hyp3} (with \eqref{PStr2a}),
and \eqref{PStr5c}, we have  
\begin{align}
\Big\| \|\jb{\nb}^{-\eps} :\! z_N^\l\!: \|_{L^q_TL^r_x}\bigg\|_{L^p(\O)}
& \leq \Big\| \|\jb{\nb}^{-\eps}:\! z_N^\l(t,x)\!: \|_{L^p(\O)}\bigg\|_{L^q_TL^r_x} \notag \\
&  \leq  C_\l\,  p^\frac{\l}{2}
\Big\| \|\jb{\nb}^{-\eps} :\! z_N^\l(t,x)\!: \|_{L^2(\O)}\Big\|_{L^q_TL^r_x}  \notag \\
& \les T^\frac{1}{q} p^\frac{\l}{2}, 
\label{PStr6}
\end{align}

\noi
for all $p \geq \max(q, r)$.
Finally,  \eqref{PStr3} follows from
\eqref{PStr6} and  
Chebyshev's inequality.

A similar computation with  Lemma \ref{LEM:W1}, \eqref{W4}, and Lemma \ref{LEM:hyp3}
shows that the sequence $\big\{:\! z_N^\l\!:\big\}_{N \in \N}$ is a Cauchy sequence  
in $L^p(\O; L^q([-T, T]; W^{-\eps, r}(\T^2)) )$.
\end{proof}

\begin{remark}\label{REM:conti}\rm
As a corollary to Proposition~\ref{PROP:PStr}, 
we can show that the tail estimate \eqref{PStr3}
and the convergence of 
$:\! z_N^\l\!:$ to $:\! z^\l\!:$
hold even when $q = \infty$ and/or $r = \infty$.
This follows from applying Sobolev's inequality (in time and/or space)
and using the fact that 
$z$ solves the linear wave/Klein-Gordon equation.  See \cite{BTT0}.
With this observation, we can easily show  that 
$:\! z_N^\l\!:, \ :\! z^\l\!: \, 
\in C([-T, T]; W^{-\eps, r}(\T^2)) $ almost surely.
See also \cite{GKO, OP}.
\end{remark}

\section{Local well-posedness of the Wick ordered NLW}
\label{SEC:LWP}

In this section, we present the proof of Theorem~\ref{THM:LWP}.
We combine the deterministic analysis via the Fourier restriction norm method
(with the hyperbolic Sobolev spaces)
and the stochastic estimate on the Wick ordered monomials 
$:\! z^\l\!:$ 
(Proposition~\ref{PROP:PStr}).
In the following, we fix $\rho > 0$.

\subsection{Hyperbolic Sobolev spaces and  Strichartz estimates}

We first  recall the hyperbolic Sobolev space $X^{s, b}$
due to Klainerman-Machedon \cite{KM} and Bourgain \cite{BO93}, 
defined by the norm
\begin{align*}
\|u\|_{X^{s, b} (\R \times \T^2)} = \|\jb{n}^s \jb{|\tau|- \jb{n}_\rho}^b \ft u( \tau,n)\|_{\l^2_n L^2_\tau( \R \times \Z^2)}.
\end{align*}

\noi
For $ b > \frac{1}{2}$, we have $X^{s, b} \subset C(\R; H^s)$.
Given an  interval $I \subset \R$,
we define the local-in-time version $X^{s, b}(I)$
as a restriction norm:
\[ \|u \|_{X^{s, b}(I)} = \inf\big\{ \|v\|_{X^{s, b}(\R \times \T^2 )}: \, v|_I = u\big\}.\]

\noi
When $I = [-T, T]$, we set $X^{s, b}_T = X^{s, b}(I)$.

The main deterministic tool for the proof of Theorem~\ref{THM:LWP}
is the following Strichartz estimates
for the  linear wave/Klein-Gordon  equation. 
Given  $0 \leq s \leq 1$, 
we say that a pair $(q, r)$ is $s$-admissible
if $2 < q \leq \infty$, 
 $2 \leq r < \infty$, 
\begin{align*}
\frac{1}{q} + \frac 2r  =  1 - s, 
\qquad \text{and} \qquad 
\frac 1q + \frac{1}{2r} \leq \frac 1 4.
\end{align*}

\noi
Then, we have the following Strichartz estimates.

\begin{lemma}\label{LEM:Str1}
Let $T \leq 1$. Given $0 \leq s \leq 1$,
let $(q, r)$ be $s$-admissible.
Then, we have 
\begin{align}
\|S(t)(\phi_0, \phi_1)  \|_{L^q_TL^r_x(\T^2)}
\les
\|(\phi_0 ,\phi_1) \|_{\H^{s}(\T^2)}.
\label{Str1}
\end{align}

\end{lemma}

See Ginibre-Velo \cite{GV}, Lindblad-Sogge \cite{LS},  and Keel-Tao \cite{KeelTao}
for the Strichartz estimates  on $\R^d$.
See also \cite{KSV}. 
The Strichartz estimates \eqref{Str1} on $\T^2$ in Lemma \ref{LEM:Str1} follows from 
those on $\R^2$ and the finite speed of propagation.

When $b > \frac 12$, the $X^{s, b}$-spaces enjoy the transference principle.
In particular, 
as a corollary to Lemma \ref{LEM:Str1}, 
we obtain the following space-time estimate.
See \cite{KS, TAO} for the proof.

\begin{lemma}\label{LEM:Str2}
Let $T \leq 1$. 
Given $0 \leq s \leq 1$,
let $(q, r)$ be $s$-admissible.
Then, for $b > \frac{1}{2}$, we have 
\begin{align*}
\| u \|_{L^q_T L^r_x} \les
\| u\|_{X_T^{s, b}}.
\end{align*}

\end{lemma}

Lastly, we state the nonhomogeneous linear  estimate.
See  \cite{GTV}.

\begin{lemma}\label{LEM:lin1}
Let $ - \frac 12 < b' \leq 0 \leq b \leq b'+1$.
Then, for $T \leq 1$, we have 
\begin{align*}
\bigg\|  \int_0^t \frac{\sin((t-t')\jb{\nb}_\rho)}{\jb{\nb}_\rho} F(t') dt'\bigg\|_{X^{s, b}_T}
\les T^{1-b+b'} \|F\|_{X^{s-1, b'}_T}.
\end{align*}

\end{lemma}

\subsection{Proof of Theorem~\ref{THM:LWP}}

In the following, we simply consider the case  $s =  s_\text{crit}+\delta$ with $\delta \ll 1$.
Given $T \leq 1$, define $\Psi(w)$
by 
\[ \Psi(w)(t) = \Psi^\o(w)(t)
 \stackrel{\text{def}}{=} 
 \int_0^t \frac{\sin ((t - t') \jb{\nb}_\rho)}{\jb{\nb}_\rho} 
  :\! (w+z)^{2m+1}(t')\!:  dt'.
\]

\noi
Let $b =  \frac 12+$. Then, for $0< \theta \leq 1-b$, 
by Lemma \ref{LEM:lin1}, we have 
\begin{align}
\| \Psi(w) \|_{X^{s, b}_T } 
\les
 T^\theta \| :\! (w+z)^{2m+1}\!\!:  \|_{X^{s-1, b-1+\theta}_T}.
\label{LWP1}
\end{align}

\noi
From \eqref{Herm4}, we have 
\begin{align*}
 :\! (w+z)^{2m+1}\!\!:  \, \, 
  = \sum_{\l = 0}^{2m+1}
\begin{pmatrix}
2m +1 \\ \l
\end{pmatrix}
 w^{2m+1 - \l}  :\! z^\l\!:.
\end{align*}

\noi
Then, 
by duality, we have 
\begin{align}
\| :\! (w+z)^{2m+1} & \!\!:  \|_{X_T^{s-1, b-1+\theta}}
 \leq 
\sum_{\l = 0}^{2m+1}C_{m, \l}
\| \, w^{2m+1-\l} :\! z^\l\!: \|_{X_T^{s-1, b-1+\theta}}\notag \\
& \leq 
\sum_{\l = 0}^{2m+1}C_{m, \l}
\sup_{h_\l}
\bigg|\iint  \ind_{[-T, T]} \, \wt w^{2m+1-\l} :\! z^\l\!:  h_\l \, dx dt\bigg|, 
\label{LWP2}
\end{align}

\noi
for any extension $\wt w$ of $w$, 
where the supremum is taken over $h_\l \in X^{1-s, 1 - b - \theta}$ 
with $\|h_\l\|_{X^{1-s, 1 - b - \theta}} = 1$.
By choosing $\theta > 0$ sufficiently small, 
we have $1 - b - \theta = \frac 12 -$.

\medskip

\noi
$\bullet$ {\bf Case 1:} $m = 1$.
\\
\indent
In this case, we have 
$s = s_\text{crit} + \dl = \frac 14 + \dl$.
Noting that  $ \big(\frac{12}{1+2\dl}, \frac{3}{1-\dl}\big)$
is $\big(\frac{1}{4} + \frac{1}{2}\dl\big)$-admissible, 
it follows from Lemma \ref{LEM:Str2} that 
\begin{align*}
\|\jb{\nb}^\eps \wt w\|_{L_T^\frac{12}{1+2\dl} L^\frac{3}{1-\dl}_x}
& \les \|\wt w\|_{X^{\frac 14 + \frac 12 \dl + \eps, \frac 12 +}}
 \les \|\wt w\|_{X^{s, \frac 12 +}}
\end{align*}

\noi
for any extension $\wt w$ of $w$, 
as long as $\eps \leq \frac 12 \dl$. 	
By taking an infimum over all the extensions $\wt w$ of $w$, we obtain
\begin{align}
\inf_{\wt w|_{[-T, T]} = w} \|\jb{\nb}^\eps  \wt w\|_{L_T^\frac{12}{1+2\dl} L^\frac{3}{1-\dl}_x}
 \les \| w\|_{X_T^{s, \frac 12 +}}.
\label{LWP3}
\end{align}

On the one hand, noting that 
 $\big(\frac{4}{1-2\dl}, \frac{1}{\dl}\big)$
is $\big(\frac{3}{4} - \frac{3}{2}\dl\big)$-admissible,
H\"older's inequality (with $T\leq 1$) and 
Lemma \ref{LEM:Str2} yield
\begin{align}
\| \jb{\nb}^\eps h_\l\|_{L^{\frac{4}{3-2\dl}}_TL^\frac{1}{\dl}_x}
& 
\les \| \jb{\nb}^\eps h_\l\|_{L^{\frac{4}{1-2\dl}}_TL^\frac{1}{\dl}_x}
\les\| \jb{\nb}^\eps h_\l\|_{X^{\frac 34 -\frac 32 \dl, \frac 12 + }}.
\label{LWP4}
\end{align}

\noi
On the other hand, applying H\"older's inequality in $t$ and   Sobolev's inequality in $x$, we have 
\begin{align}
\| \jb{\nb}^\eps h_\l\|_{L^{\frac{4}{3-2\dl}}_TL^\frac{1}{\dl}_x}
\les \| \jb{\nb}^\eps h_\l\|_{L^{2}_TL^\frac{1}{\dl}_x}
\les \| \jb{\nb}^\eps h_\l\|_{X^{1 -2 \dl, 0}}.
\label{LWP5}
\end{align}

\noi
Interpolating \eqref{LWP4} and \eqref{LWP5}
with  sufficiently small $\theta > 0$, we obtain
\begin{align}
\| \jb{\nb}^\eps h_\l\|_{L^{\frac{4}{3-2\dl}}_TL^\frac{1}{\dl}_x}
& \les \|  h_\l\|_{X^{\frac{3}{4}  -\frac{5}{4} \dl + \eps, 1-b-\theta}}
 \les \|  h_\l\|_{X^{1-s, 1-b-\theta}}
\label{LWP6}
\end{align}

\noi
as long as $\eps \leq \frac 14 \dl$.

For $\l = 0, 1, 2, 3$, define $(q_\l, r_\l)$
by 
\[ 1 = (3-\l) \frac{1+2\dl}{12} +\frac{3-2\dl}{4} + \frac{1}{q_\l}
\qquad \text{and}\qquad 
1 = (3-\l) \frac{1-\dl}{3} + \dl + \frac{1}{r_\l}.\]

\noi
When $\l = 0$, we have  $q_0 = r_0 = \infty$ and $:\! z^0\!: \ \equiv 1$.
Then, by fractional Leibniz rule and H\"older's inequality with \eqref{LWP3} and \eqref{LWP6}, 
we have 
\begin{align}
\inf_{\wt w|_{[-T, T]} = w}\bigg|\iint   
& \ind_{[-T, T]} \wt w^{3-\l}  :\! z^\l\!: h_\l \, dx dt\bigg| \notag\\
& = \inf_{\wt w|_{[-T, T]} = w}\bigg|\iint   
\ind_{[-T, T]} \jb{\nb}^\eps (\wt w^{3-\l}h_\l) \jb{\nb}^{-\eps} :\! z^\l\!: dx dt \bigg|\notag \\
& \les \inf_{\wt w|_{[-T, T]} = w} \|\jb{\nb}^\eps \wt w\|_{L^\frac{12}{1+2\dl}_T L^\frac{3}{1-\dl}_x}^{3-\l}
\| \jb{\nb}^\eps h_\l\|_{L^{\frac{4}{3-2\dl}}_TL^\frac{1}{\dl}_x}
\|\jb{\nb}^{-\eps} :\! z^\l\!: \|_{L^{q_\l}_TL^{r_\l}_x}\notag \\
& \les  \|w\|_{X^{s, b}_T}^{3-\l} \|  h_\l\|_{X^{1-s, 1-b-\theta}}
\|\jb{\nb}^{-\eps} :\! z^\l\!: \|_{L^{q_\l}_TL^{r_\l}_x}\notag\\
& =  \|w\|_{X^{s, b}_T}^{3-\l}
\|\jb{\nb}^{-\eps} :\! z^\l\!: \|_{L^{q_\l}_TL^{r_\l}_x}
\label{LWP7}
\end{align}

\noi
as long as $0 < \eps \leq \frac 14 \dl$.
Hence, by Proposition~\ref{PROP:PStr} with \eqref{LWP1}, \eqref{LWP2},  and \eqref{LWP7}, we obtain
\begin{align}
\|\Psi(w)\|_{X^{s, b}_T}
 \les T^\theta
\sum_{\l = 0}^{3}
\|w\|_{X^{s, b}_T}^{3-\l}
\label{LWP8}
\end{align}

\noi
outside a set of probability $< \exp\big(-\frac{1}{T^c}\big)$ for some $c > 0$.
Similarly, we have 
\begin{align}
\|\Psi(w_1) - \Psi(w_2) \|_{X^{s, b}_T}
 \les T^\theta
\sum_{\l = 0}^{2}
\big(\|w_1\|_{X^{s, b}_T}^{2-\l} + \|w_2\|_{X^{s, b}_T}^{2-\l} \big)\|w_1 -w_2\|_{X^{s, b}_T}
\label{LWP9}
\end{align}

\noi
outside a set of probability $< \exp\big(-\frac{1}{T^c}\big)$.
Therefore, it follows from \eqref{LWP8} and \eqref{LWP9}
that 
for each $T\ll1$, 
there exists
 a set $\O_T$ with $P(\O_T^c) < \exp\big(-\frac{1}{T^c}\big)$
such that, for each $\o \in \O_T$,  
 $\Psi^\o$ is a 
contraction on a ball of radius $O(1)$ in $X^{s, b}_T$.

\medskip

\noi
$\bullet$ {\bf Case 2:} $m \geq 2$.
\\
\indent
In this case, we have 
$s = s_\text{crit} + \dl = 1-\frac 1m  + \dl$.
Define $(q, r)$ by 
\[ \frac 1q = \frac{3m-1}{3m(2m+1)} + \frac{\dl}{6} 
\qquad 
\text{and} \qquad
 \frac 1r = \frac{3m+4}{6m(2m+1)} - \frac{\dl}{3}.
\]

\noi
Noting that $(q, r)$ 
is $\big(s_\text{crit} + \frac{1}{2}\dl\big)$-admissible, 
it follows from Lemma \ref{LEM:Str2} that 
\begin{align*}
\|\jb{\nb}^\eps \wt w\|_{L^q_T L^r_x}
& \les \|\wt w\|_{X^{s_\text{crit} + \frac 12 \dl + \eps, \frac 12 +}}
 \les \|\wt w\|_{X^{s, \frac 12 +}}
\end{align*}

\noi
for any extension $\wt w$ of $w$, 
as long as $\eps \leq \frac 12 \dl$. 	
By taking an infimum over all the extensions $\wt w$ of $w$, we obtain
\begin{align}
\inf_{\wt w|_{[-T, T]} = w} \|\jb{\nb}^\eps \wt w\|_{L^q_T L^r_x}
&  \les \|w\|_{X^{s, \frac 12 +}_T}.
\label{LWP10}
\end{align}

Now, define $(\wt q,\wt r)$ by 
\[ \frac 1{\wt q} = \frac{1}{3m} - \frac{2m+1}{6}\dl 
\qquad 
\text{and} \qquad
 \frac 1{\wt r} = \frac{3m-4}{6m} + \frac{2m+1}{3}\dl.
\]

\noi
Then,  $(\wt q, \wt r)$ is  $\big(1-s_\text{crit}   - \frac{2m+1}{2}\dl\big)$-admissible.
On the one hand, 
by Lemma \ref{LEM:Str2}, we have 
\begin{align}
\| \jb{\nb}^\eps h_\l\|_{L^{\wt q}_TL^{\wt r}_x}
& 
\les\| \jb{\nb}^\eps h_\l\|_{X^{1-s_\text{crit} - \frac{2m+1}{2}\dl, \frac 12 + }}.
\label{LWP11}
\end{align}

\noi
On the other hand,  by  Sobolev's inequality, we have 
\begin{align}
\| \jb{\nb}^\eps h_\l\|_{L^{\wt q}_TL^{\wt r}_x}
\les \| \jb{\nb}^\eps h_\l\|_{X^{1 -\frac{3m-4}{3m} - \frac{4m+2}{3} \dl, \frac{1}{2}-\frac{1}{3m}+\frac{2m+1}{6}\dl}}.
\label{LWP12}
\end{align}

\noi
Note that the temporal regularity on the right-hand side of \eqref{LWP12}
is less than $\frac 12$ by choosing $\dl > 0$ sufficiently small.
Hence, by interpolating \eqref{LWP11} and \eqref{LWP12}
with sufficiently small $\theta > 0$, we obtain
\begin{align}
\| \jb{\nb}^\eps h_\l\|_{L^{\wt q}_TL^{\wt r}_x}
& \les \| h_\l\|_{X^{1-s_\text{crit}  -m \dl +\eps, 1-b-\theta}}
 \les \|  h_\l\|_{X^{1-s, 1-b-\theta}}
\label{LWP13}
\end{align}

\noi
as long as $\eps \leq  (m-1) \dl$.

Proceeding as before, it follows from  H\"older's inequality with \eqref{LWP10} and \eqref{LWP13}
that 
\begin{align}
\bigg|\iint   
\inf_{\wt w|_{[-T, T]} = w} & \ind_{[-T, T]} \wt w^{2m+1-\l}  :\! z^\l\!:  h_\l \, dx dt\bigg| \notag\\
& = \inf_{\wt w|_{[-T, T]} = w} \bigg|\iint   
\ind_{[-T, T]} \jb{\nb}^\eps (\wt w^{2m+1-\l}h_\l) \jb{\nb}^{-\eps} :\! z^\l\!: dx dt\bigg|\notag \\
& \les\inf_{\wt w|_{[-T, T]} = w} \|\jb{\nb}^\eps \wt w\|_{L^q_T L^r_x}^{2m+1-\l}
\| \jb{\nb}^\eps h_\l\|_{L^{\wt q}_TL^{\wt r}_x}
\|\jb{\nb}^{-\eps} :\! z^\l\!: \|_{L^{\frac q\l}_TL^{\frac r\l}_x}\notag \\
& \les  \|w\|_{X^{s, b}_T}^{2m+1-\l} \|  h_\l\|_{X^{1-s, 1-b-\theta}}
\|\jb{\nb}^{-\eps} :\! z^\l\!: \|_{L^{\frac q\l}_TL^{\frac r\l}_x}\notag\\
& =  \|w\|_{X^{s, b}_T}^{2m+1-\l}
\|\jb{\nb}^{-\eps} :\! z^\l\!: \|_{L^{\frac q\l}_TL^{\frac r\l}_x}
\label{LWP14}
\end{align}

\noi
as long as $0 < \eps \leq \frac 12 \dl$.
Hence, by Proposition~\ref{PROP:PStr} with \eqref{LWP1}, \eqref{LWP2},  and \eqref{LWP14}, we obtain
\begin{align*}
\|\Psi(w)\|_{X^{s, b}_T}
&  \les T^\theta
\sum_{\l = 0}^{2m+1}
\|w\|_{X^{s, b}_T}^{2m+1-\l}, 
\\
\|\Psi(w_1) - \Psi(w_2) \|_{X^{s, b}_T}
&  \les T^\theta
\sum_{\l = 0}^{2m}
\big(\|w_1\|_{X^{s, b}_T}^{2m-\l} + \|w_2\|_{X^{s, b}_T}^{2m-\l} \big)\|w_1 - w_2\|_{X^{s, b}_T}
\end{align*}

\noi
outside a set of probability $< \exp\big(-\frac{1}{T^c}\big)$ for some $c > 0$.
Therefore, 
for each $T\ll1$, 
there exists
 a set $\O_T$ with $P(\O_T^c) < \exp\big(-\frac{1}{T^c}\big)$
such that, for each $\o \in \O_T$,  
 $\Psi^\o$ is a 
contraction on a ball of radius $O(1)$ in $X^{s, b}_T$.

This completes the proof of Theorem~\ref{THM:LWP}.

\section{Weak universality: Wick ordered NLW as a scaling limit}

In this section, we present the proof of Theorem~\ref{THM:WNLW}.
We follow closely the argument in~\cite{GKO}.
With $z_\eps  = z_\eps^\o = S(t) (\phi^\o_{\eps, 0}, \phi^\o_{\eps, 1})$, 
let us decompose $u_{\eps} = z_{\eps} +w_{\eps}$ 
as in \eqref{decomp2}.
Then, the residual term $w_{\eps}$ satisfies 
\begin{align}
 \dt^2 w_\eps - \Dl w_\eps + \rho w_\eps = F_\eps(w_\eps), 
\label{XNLW3}
 \end{align}

\noi
where $F_\eps(w_\eps)$ is given by 
\begin{align}
 F_{\eps} (w_{\eps}) & = \eps^{-3} 
\big\{f (\eps (z_{\eps} +   w_{\eps})) + 
\eps (\eps^2 \rho -  \rho_\eps) (z_\eps + w_\eps) \big\} \notag\\
& = \eps^{- 2} \{f' (0) + \eps^2 \rho -  \rho_\eps\} (z_{\eps} + w_{\eps}) 
+ \frac{f^{(3)} (0)}{6}  (z_{\eps} +w_{\eps})^3 + R_{\eps}, 
\label{MG1}
\end{align}

\noi
where the second equality follows 
from $f (0) = f'' (0) = 0$
and  Taylor's remainder theorem 
with the remainder term $R_\eps$ given by 
\begin{align}
 R_{\eps} =  \eps \int_0^1 \frac{(1 -
   \theta)^3}{6} f^{(4)} (\theta \eps (z_{\eps} +   w_{\eps})) d \theta 
   \cdot (z_{\eps} +   w_{\eps})^4 . 
\label{MG2}
\end{align}

From \eqref{A2}, we see that $z_\eps(t,x)$ 
 is a mean-zero real-valued Gaussian random variable with variance
 \[\s_{\eps} =  \E[z_{\eps}^2 (t,x)] \sim  \log \eps^{-1} .\]

\noi
Note that $\s_\eps$ is independent of $x \in \T^2$ and $t\in \R$.
Recalling from \eqref{H1a} that $x^3 = H_3(x; \s) + 3\s x$, 
it follows from \eqref{MG1} and \eqref{MG2} that
\begin{align*}
 F_{\eps} (w_{\eps}) = \eps^{- 2} \bigg\{ f' (0) +  \eps^2 \rho -  \rho_\eps
 & + 3 \eps^{2 } \sigma_{\eps}
\frac{f^{(3)} (0)}{6} \bigg\} (z_{\eps} + w_{\eps}) \notag\\ 
& +  \frac{f^{(3)} (0)}{6} H_3 (z_{\eps} +w_{\eps}; \sigma_{\eps}) + R_{\eps}. 
\end{align*}

\noi
For each $\eps > 0$, we set $\rho_\eps$ by 
\[ \rho_\eps =   f' (0) + \eps^2 \rho  + \eps^2   \s_{\eps}  \frac{f^{(3)} (0)}{2} \]

\noi
so that the first term on the right-hand side vanishes.
Then, by letting 
$\ld =\frac{f^{(3)} (0)}{6} $,  we obtain 
\begin{align*}
 F_{\eps} (w_{\eps}) 
= \lambda H_3 (z_{\eps} +  w_{\eps}; \sigma_{\eps}) + R_\eps
 \stackrel{\text{def}}= \ld  :\! u_\eps^3\!: + R_\eps. 
\end{align*}

From  \eqref{MG2}
and \eqref{bound}, we have
\begin{align*}
| R_{\eps} |
& = \bigg|\eps  \int_0^1 \frac{(1 -   \theta)^3}{6} 
 f^{(4)} ( \theta \eps (z_{\eps} +   w_{\eps}) )
 d \theta    \cdot (z_{\eps} +   w_{\eps})^4 \bigg|\notag\\
& \les \eps \big\{|z_{\eps}| +   |w_{\eps}|\big\}^{M+4}.
\end{align*}

\noi
In particular, we can write \eqref{XNLW3} as 
\begin{align}
 \dt^2 w_\eps - \Dl w_\eps  + \rho w_\eps 
& = \ld
  \sum_{\l=0}^3 {3\choose \l} :\! z_\eps^\l \,  \!: w_\eps^{3-\ell}
  + O\big(  \eps \big\{|z_{\eps}| +   |w_{\eps}|\big\}^{M+4}\big).
\label{XNLW4}
 \end{align}

\noi
It follows from 
 Proposition~\ref{PROP:PStr} with \eqref{A2}
that 
\begin{align*}
 \eps \| z_{\eps} \|_{L^{q}_T L^r_x}^{M+4} =
 o_{\eps} (1) 
\end{align*}

\noi
almost surely.
Then, by proceeding  as in Section \ref{SEC:LWP}
(
where we handle the second term on the right-hand side
of \eqref{XNLW4} by applying the argument in 
Section \ref{SEC:LWP} with $2m+1 \geq M+4$), 
 we obtain an a priori bound on $w_\eps$, uniformly in $\eps > 0$.
Moreover, the local existence time $T = T_\o $ 
can be chosen to be  independent of $\eps>0$.

Let $u$ be the solution to \eqref{NLW_ld}.
In an analogous manner, we can estimate the difference $w -w_\eps$, 
where $w = u - z$ as in \eqref{decomp2}.
Together with the almost sure convergence of~$z_\eps$ to~$z$
(see Remark \ref{REM:conti}), 
we see that $u_\eps$ converges to $u$ in $C([-T_\o, T_\o]; H^s(\T^2))$
for $s < 0$.
This completes the proof of Theorem~\ref{THM:WNLW}.

\begin{ackno}\rm
T.O.~is supported by the ERC starting grant 
no.~637995 ``ProbDynDispEq''.
L.T.~is supported by the grants "BEKAM"  ANR-15-CE40-0001 and "ISDEEC'' ANR-16-CE40-0013. 

\end{ackno}


\begin{thebibliography}{99}

\bibitem{AC}
S.~Albeverio, A.~Cruzeiro,
{\it  Global flows with invariant (Gibbs) measures for Euler and Navier-Stokes two dimensional fluids,} 
Comm. Math. Phys. 129 (1990) 431--444.




\bibitem{BO93}
J.~Bourgain, 
{\it Fourier transform restriction phenomena for certain lattice subsets and applications to nonlinear evolution equations, I: Schr\"odinger equations,} Geom. Funct. Anal. 3 (1993), 107--156.


\bibitem{BO94}
J.~Bourgain, 
{\it Periodic nonlinear Schr\"odinger equation and invariant measures}, 
Comm. Math. Phys. 166 (1994), no. 1, 1--26.

\bibitem{BO96}
J.~Bourgain, 
{\it Invariant measures for the 2D-defocusing nonlinear Schr\"odinger equation}, 
Comm. Math. Phys. 176 (1996), no. 2, 421--445. 

\bibitem{BoPCMI}
J.~Bourgain, 
{\it Nonlinear Schr\"odinger equations,}
 Hyperbolic equations and frequency interactions (Park City, UT, 1995), 3--157, IAS/Park City Math. Ser., 5, Amer. Math. Soc., Providence, RI, 1999.


\bibitem{BS}
D.~Brydges, G.~Slade, 
{\it Statistical mechanics of the 2-dimensional focusing nonlinear Schr\"odinger equation,}
 Comm. Math. Phys. 182 (1996), no. 2, 485--504.

\bibitem{BTT0}
N.~Burq, L.~Thomann, N.~Tzvetkov, 
{\it Global infinite energy solutions for the cubic wave equation,}
 Bull. Soc. Math. France 143 (2015), no. 2, 301--313. 

\bibitem{BTT1}
N.~Burq, L.~Thomann, N.~Tzvetkov, 
{\it Remarks on the Gibbs measures for nonlinear dispersive equations}, 
 Ann. Fac. Sci. Toulouse Math., to appear.


\bibitem{BT1}
N.~Burq, N.~Tzvetkov,
{\it Random data Cauchy theory for supercritical wave equations. I. Local theory,}
Invent. Math. 173 (2008), no. 3, 449--475.


\bibitem{BT2}
N.~Burq, N.~Tzvetkov, 
{\it Random data Cauchy theory for supercritical wave equations. II. A global existence result,} Invent. Math. 173 (2008), no. 3, 477--496.




\bibitem{DPD}
G.~Da Prato, A.~Debussche, 
{\it Two-dimensional Navier-Stokes equations driven by a space-time white noise}, 
J. Funct. Anal. 196 (2002), no. 1, 180--210.

\bibitem{DPT1}
G.~Da Prato, L.~Tubaro, 
{\it Wick powers in stochastic PDEs: an introduction,} 
Technical Report UTM, 2006, 39 pp.

\bibitem{GTV}
J.~Ginibre, Y.~Tsutsumi, G.~Velo, 
{\it On the Cauchy problem for the Zakharov system,} J. Funct. Anal. 151 (1997), no. 2, 384--436.

\bibitem{GV}
J.~Ginibre, G.~Velo, 
{\it Generalized Strichartz inequalities for the wave equation,} J. Funct. Anal. 133 (1995),
50--68.

\bibitem{GJ}
J.~Glimm, A.~Jaffe, 
{\it Quantum physics. A functional integral point of view,} Second edition. Springer-Verlag, New York, 1987. xxii+535 pp.



\bibitem{GKO}
M.~Gubinelli, H.~Koch, T.~Oh,
{\it  Renormalization of the two-dimensional stochastic nonlinear wave equations}, arXiv:1703.05461  [math.AP].


\bibitem{GP1}
M.~Gubinelli, N.~Perkowski, 
{\it   KPZ reloaded,},
Commun. Math. Phys. 349 (2017), 165--269.


\bibitem{GP2}
M.~Gubinelli, N.~Perkowski, 
{\it The Hairer--Quastel universality result at stationarity}, 
 to appear in RIMS K\^oky\^uroku Bessatsu (2016).


\bibitem{HQ} 
M.~Hairer, J.~Quastel,
{\it A class of growth models rescaling to KPZ,}
arXiv:1512.07845 [math-ph].



\bibitem{KeelTao}
M.~Keel, T.~Tao, {\it Endpoint Strichartz estimates},  Amer. J. Math. 120 (1998), no. 5, 955--980.




\bibitem{KSV}
R.~Killip, B.~Stovall, M.~Vi\c{s}an, 
{\it Blowup behaviour for the nonlinear Klein-Gordon equation}, 
Math. Ann. 358 (2014), no. 1-2, 289--350.




\bibitem{KM}
S.~Klainerman, M.~Machedon, 
{Space-time estimates for null forms and the local existence theorem}, Comm. Pure Appl. Math., 46 (1993), 1221--1268.



\bibitem{KS}
S.~Klainerman, S.~Selberg, {\it Bilinear estimates and applications to nonlinear wave equations,} Commun. Contemp. Math. 4 (2002), no. 2, 223--295.


\bibitem{LS}
H.~Lindblad, C.~Sogge, {\it On existence and scattering with minimal regularity for semilinear wave equations,} J. Funct. Anal. 130 (1995), 357--426.


\bibitem{McKean}
H.P.~McKean, 
{\it Statistical mechanics of nonlinear wave equations. IV. Cubic Schr\"odinger,} 
 Comm. Math. Phys. 168 (1995), no. 3, 479--491. 
 {\it Erratum: Statistical mechanics of nonlinear wave equations. IV. Cubic Schr\"odinger}, Comm. Math. Phys. 173 (1995), no. 3, 675.


\bibitem{Nelson2}
E.~Nelson, 
{\it A quartic interaction in two dimensions}, 
 1966 Mathematical Theory of Elementary Particles (Proc. Conf., Dedham, Mass., 1965) pp. 69--73 M.I.T. Press, Cambridge, Mass.






\bibitem{OP}
T.~Oh, O.~Pocovnicu, 
{\it Probabilistic global well-posedness of the energy-critical defocusing quintic nonlinear wave equation 
on $\R^3$}, J. Math. Pures Appl. 105 (2016), 342--366. 



\bibitem{OQ}
T.~Oh, J.~Quastel, 
{\it On Cameron-Martin theorem and almost sure global existence}, 
Proc. Edinb. Math. Soc. 59 (2016), 483--501. 


\bibitem{ORT}
T.~Oh, G.~Richards, L.~Thomann, 
{\it On invariant Gibbs measures for the generalized KdV equations}, 
Dyn. Partial Differ. Equ. 13 (2016), no.2, 133--153. 

\bibitem{OT}
T.~Oh, L.~Thomann, 
{\it 
 A pedestrian approach to the invariant Gibbs measure for 
the 2-$d$ defocusing nonlinear Schr\"odinger equations}, 
arXiv:1509.02093 [math.AP].





\bibitem{Simon}
B.~Simon, 
{\it  The $P(\varphi)_2$ Euclidean (quantum) field theory,} Princeton Series in Physics. Princeton University Press, Princeton, N.J., 1974. xx+392 pp.


\bibitem{TAO}
T.~Tao, {\it Nonlinear dispersive equations. Local and global analysis,} CBMS Regional Conference Series in Mathematics, 106. Published for the Conference Board of the Mathematical Sciences, Washington, DC; by the American Mathematical Society, Providence, RI, 2006. xvi+373 pp.







\bibitem{Zworski}
M.~Zworski,    {\it Semiclassical analysis.} Graduate Studies in Mathematics, 138. American Mathematical Society, Providence, RI, 2012. xii+431 pp.


\end{thebibliography}
\end{document}